\documentclass[12pt,a4paper]{article}

\usepackage{amsmath}
\usepackage{amssymb}
\usepackage{amsthm}

\numberwithin{equation}{section}

\begin{document}

\title{A coloring of the square of the $8$-cube with $13$ colors}

\author{%
Janne I. Kokkala\footnote{Supported 
by the Aalto ELEC Doctoral School, by the Nokia Foundation, and by the Emil Aaltonen Foundation.}\ \ 
and Patric R. J. \"Osterg\aa rd\\
\hspace*{5mm}\\
Department of Communications and Networking\\
Aalto University School of Electrical Engineering\\
P.O.\ Box 13000, 00076 Aalto, Finland
}

\date{}

\maketitle

\begin{abstract}
Let $\chi_{\bar{k}}(n)$ be the number of colors required to color the $n$-dimensional hypercube such that no two vertices with the same color are at a distance at most $k$. In other words, $\chi_{\bar{k}}(n)$ is the minimum number of binary codes with minimum distance at least $k+1$ required to partition the $n$-dimensional Hamming space. By giving an explicit coloring, it is shown that $\chi_{\bar{2}}(8)=13$.
\end{abstract}

\section{Introduction}

For any pair $u, v \in \{0,1\}^n$, the \emph{Hamming distance} between $u$ and $v$, denoted by $d_H(u,v)$,
is the number of positions in which $u$ and $v$ differ. A \emph{binary} $(n,M,d)$ \emph{code} $C$ is a subset of $\{0,1\}^n$ for which $|C|=M$ and the minimum Hamming distance between any two distinct elements of $C$ is $d$. The parameters $n$, $M$, and $d$ are called the \emph{length}, the \emph{size}, and the \emph{minimum distance} of $C$, respectively.

The \emph{$n$-dimensional hypercube}, also called the \emph{$n$-cube}, denoted by $Q_n$, is the graph with vertex set $V=\{0,1\}^n$ such that two vertices are adjacent if and only if their Hamming distance is exactly $1$. Given a graph $G$, the \emph{$k$th power} of $G$, denoted by $G^k$, is the graph obtained from $G$ by adding edges between all pairs of vertices that have distance at most $k$ in $G$. In particular, $G^2$ is called the \emph{square} of $G$.

A proper vertex coloring of $Q_n^k$ corresponds to a partition of $\{0,1\}^n$ into binary codes of minimum distance at least $k+1$. The chromatic number of $Q_n^k$ is denoted by $\chi_{\bar{k}}(n)$. The problem of finding bounds and exact values for $\chi_{\bar{k}}(n)$ arises from the problem of scalability of certain optical networks and has attracted wide interest in coding theory and combinatorics; see for example \cite{W97, KDP00, Z01, NDG02, O04}.

\section{Determining  $\chi_{\bar{2}}(8)$}

The size of a binary code of length $8$ and minimum distance $3$ is at most $20$ \cite{BBMOS78}. Therefore, at least $\left\lceil \frac{2^8}{20} \right\rceil = 13$ colors are needed to color the square of the $8$-cube. Colorings with $14$ colors were first obtained by Hougardy in 1991 \cite{Z01} and Royle in 1993 \cite[Section~9.7]{JT95}, but it has been an open problem whether $13$ colors suffice.

We give a partition of $\{0,1\}^8$ into $13$ codes of minimum distance at least $3$ in Table~\ref{tab:coloring}, which shows that $\chi_{\bar{2}}(8)=13$. To save space, the elements of $\{0,1\}^8$ are given as integers from $0$ to $255$. Twelve of the codes are $(8,20,3)$ codes and the remaining code is an $(8,16,4)$ code.

\begin{table}[t]
  \centering
  \caption{A partition of $\{0,1\}^8$ into $13$ codes of minimum distance at least $3$}
  \label{tab:coloring}
  \begin{tabular}{|r|r|r|r|r|r|r|r|r|r|r|r|r|}
\hline
$C_1$ & $C_2$ & $C_3$ & $C_4$ & $C_5$ & $C_6$ & $C_7$ & $C_8$ & $C_9$ & $C_{10}$ & $C_{11}$ & $C_{12}$ & $C_{13}$ \\
\hline
9 & 3 & 23 & 8 & 10 & 2 & 14 & 1 & 4 & 6 & 7 & 5 & 0 \\
18 & 29 & 24 & 20 & 21 & 13 & 25 & 15 & 11 & 19 & 12 & 16 & 30 \\
37 & 38 & 33 & 31 & 35 & 39 & 36 & 22 & 17 & 28 & 26 & 27 & 45 \\
56 & 40 & 47 & 46 & 44 & 52 & 55 & 42 & 54 & 32 & 41 & 34 & 51 \\
63 & 59 & 50 & 49 & 48 & 57 & 58 & 60 & 61 & 43 & 53 & 62 & 75 \\
71 & 68 & 66 & 65 & 70 & 81 & 67 & 76 & 78 & 69 & 64 & 72 & 85 \\
92 & 90 & 77 & 82 & 73 & 94 & 84 & 80 & 87 & 74 & 83 & 79 & 102 \\
96 & 109 & 100 & 103 & 95 & 107 & 104 & 91 & 88 & 89 & 93 & 86 & 120 \\
110 & 112 & 121 & 123 & 101 & 108 & 111 & 99 & 98 & 116 & 106 & 97 & 135 \\
115 & 119 & 126 & 124 & 122 & 114 & 113 & 117 & 105 & 127 & 118 & 125 & 153 \\
134 & 133 & 139 & 131 & 132 & 142 & 128 & 136 & 130 & 137 & 129 & 138 & 170 \\
149 & 150 & 140 & 154 & 143 & 144 & 141 & 157 & 158 & 159 & 148 & 156 & 180 \\
155 & 152 & 145 & 160 & 146 & 151 & 147 & 167 & 165 & 174 & 162 & 164 & 204 \\
163 & 175 & 166 & 173 & 185 & 161 & 171 & 176 & 168 & 178 & 184 & 169 & 210 \\
172 & 177 & 189 & 182 & 190 & 186 & 188 & 187 & 179 & 181 & 191 & 183 & 225 \\
202 & 201 & 199 & 196 & 209 & 197 & 198 & 194 & 193 & 192 & 203 & 195 & 255 \\
205 & 206 & 212 & 207 & 220 & 200 & 216 & 215 & 221 & 214 & 222 & 213 & \\
208 & 211 & 218 & 217 & 224 & 219 & 223 & 228 & 239 & 227 & 231 & 238 & \\
246 & 226 & 232 & 234 & 235 & 230 & 229 & 233 & 244 & 237 & 236 & 240 & \\
249 & 252 & 243 & 245 & 247 & 253 & 242 & 254 & 250 & 248 & 241 & 251 & \\
\hline
  \end{tabular}
\end{table}

The listed coloring is one of many colorings found with a computer-aided approach. The computational techniques will be discussed in detail in a full paper. It will further be checked whether these colorings can be used as substructures to obtain colorings of the square of the $9$-cube with $13$ colors.

\bibliographystyle{unsrt}
\bibliography{note}
\end{document}